\documentclass[11pt]{article}

\usepackage{fullpage}
\usepackage{multirow}
\usepackage{graphicx}
\usepackage{amsmath}
\usepackage{amssymb}
\usepackage{amsfonts}
\usepackage{natbib}
\usepackage{setspace} 
\usepackage{psfig}
\usepackage{epsfig}
\usepackage{multirow}
\usepackage{float}

\def\R{\mathbb{R}}

\begin{document} 

\title{\bf Screening and metamodeling of computer experiments with functional outputs. Application to thermal-hydraulic computations} 

\author{\bf Benjamin AUDER$^\ast$, Agn\`es DE CRECY$^\dag$, Bertrand IOOSS$^\ddag$ and Michel MARQU\`ES$^\ast$}
\date{}

\maketitle 

\begin{center}
{\it Reliability Engineering and System Safety}\\
in the special SAMO 2010 issue

\vspace{0.2cm}

  $^\ast$ CEA, DEN, Centre de Cadarache, F-13108, Saint-Paul-lez-Durance, France\\
  $^\dag$ CEA, DEN, Centre de Grenoble, F-38054, Grenoble, France\\
  $^\ddag$ EDF R\&D, 6 Quai Watier, F-78401, Chatou, France 

  \vspace{0.2cm}
  Corresponding author: B. Iooss ; Email: bertrand.iooss@edf.fr\\
  Phone: +33-1-30877969 ; Fax: +33-1-30878213
\end{center}

\doublespacing

\abstract{
To perform uncertainty, sensitivity or optimization analysis on scalar variables calculated by a cpu time expensive computer code, a widely accepted methodology consists in first identifying the most influential uncertain inputs (by screening techniques), and then in replacing the cpu time expensive model by a cpu inexpensive mathematical function, called a metamodel. This paper extends this methodology to the functional output case, for instance when the model output variables are curves. The screening approach is based on the analysis of variance and principal component analysis of output curves. The functional metamodeling consists in a curve classification step, a dimension reduction step, then a classical metamodeling step. An industrial nuclear reactor application (dealing with uncertainties in the pressurized thermal shock analysis) illustrates all these steps.
}


\section{INTRODUCTION}

The uncertainty analysis and sensitivity analysis (UASA) process is one of the key step for the development and the use of predictive complex computer models (de~Rocquigny et al. \cite{derdev08}, Saltelli et al. \cite{salcha00}). 
On one hand, this process aims at quantifying the impact of the input data and parameters uncertainties on the model predictions.
On the other hand, it investigates how the model outputs respond to variations of the model inputs. 
For example in the nuclear engineering domain, it has been applied to waste storage safety studies (Helton et al. \cite{heljoh06}), radionuclide transport modeling in the aquifer (Volkova et al. \cite{volioo08}), safety passive system evaluation (Marqu\`es et al. \cite{marpig05}) and simulation of large break loss of coolant accident scenario (de Cr\'ecy et al. \cite{decbaz08}).

In practice, when dealing with UASA methods, four main problems can arise:
\begin{enumerate}
\item	Physical models can involve complex and irregular phenomena sometimes with strong interactions between physical variables. This problem is resolved by using variance-based measures (Sobol \cite{sob93}; Saltelli et al. \cite{salcha00});
\item	Computing variance-based measures can be infeasible for cpu time consuming code. Metamodel-based techniques (Sacks et al. \cite{sacwel89}; Fang et al. \cite{fanli06}; Iooss et al. \cite{ioovan06}) solve this problem and provide a deep exploration of the model behavior.
A metamodel is a mathematical emulator function, approximating a few simulation results performed with the computer code and giving acceptable predictions with a negligible cpu time cost;
\item	The two preceeding tools (variance-based sensitivity analysis and metamodel's method) are applicable for low-dimensional problems while numerical models can take as inputs a large number of uncertain variables (typically several tens or hundreds). One simple solution is to apply, in a preliminary step, a screening technique which allows to rapidly identify the main influent input variables (Saltelli et al. \cite{salcha00}; Kleijnen \cite{kle08});
\item	Another problem of high dimensionality arises when numerical models produce functional output variables, for instance spatially or temporally dependent. Until recent years, this problem has received only little attention in the UASA framework.
However, three recent works have brought some first solutions:
\begin{itemize}
\item Campbell et al. \cite{cammck06} use a principal component analysis of output temporal curves, then compute sensitivity indices of each input on each principal component coefficient; 
 \item Lamboni et al. \cite{lammon10} develop the multivariate global sensitivity analysis method. It allows to aggregate the different sensitivity indices of the principal component coefficients in a unique index, called the generalized sensitivity index.
 Each generalized sensitivity index explains the influence of the corresponding input on the overall output curve variability.
 \item Dealing with complex and cpu time consuming computer codes, Marrel et al. \cite{marioo10} propose to build a functional metamodel (based on a wavelet decomposition technique and the Gaussian process metamodel, see Bayarri et al. \cite{bayber07}).
 Then, Monte Carlo techniques are applied on the functional metamodel to obtain variance-based sensitivity indices.
 \end{itemize}
\end{enumerate}
In this paper, we present an overall UASA methodology, from screening to metamodeling, applicable to output curves of cpu time consuming computer models. 
 
This methodology is motivated by an industrial application concerning the nuclear safety.
Such study requires the numerical simulation of the so-called pressurized thermal shock analysis using qualified computer codes. 
 This quantitative analysis aims at calculating the vessel failure probability of a nuclear pressurized reactor. 
 In a general context, structural reliability aims at determining the failure probability of a system by modeling its input variables as random variables.
The system can be simulated by a numerical model (often seen as a black box function), which can be written for instance $Y=f(X)$ where $Y \in \mathbb{R}$ is the monodimensional output variable, $f(\cdot)$ is the model function and $X=(X^1,\ldots,X^p) \in \mathbb{R}^p$ are the $p$-dimensional input variables.
A reliability problem can consist in evaluating the probability that the output $Y$ exceeds a given treshold $Y_s$.
One problem is that, as the model function (i.e. the numerical computer code) becomes more and more complex, cpu time of the model increases dramatically for each run.
To solve this problem, authors from structural reliability domain have largely contributed to the development and promotion of advanced Monte Carlo and FORM/SORM methods (Madsen et al. \cite{madkre86}).
These techniques are now widely used in other physical domains, from hydraulics to aerospace engineering (see some examples in de~Rocquigny et al. \cite{derdev08}).
For our purpose (scenario of a nuclear reactor pressurized thermal shock), past and recent works have proposed some advanced methods in order to compute the failure probability with a small number of computer code runs, additionally to a high level and conservative confidence interval (de~Rocquigny \cite {der06a}, Munoz-Zuniga et al. \cite{mungar10}).

At present, one major challenge in such calculations is to propagate input uncertainties in thermal-hydraulic models, by using for example Monte Carlo methods. 
However, the four previously described difficulties turn this problem to an extremely difficult one.
First, sophisticated phenomena are implemented in the computer code.
Then, the behaviour of the output variables are particularly complex and trivial method (as linear approximation) are not applicable.
Second, the thermal-hydraulic models require several hours to perform one simulation in order to simulate the full transient, while several thousands simulations are required by the Monte Carlo methods.
Thus, metamodels are necessary to solve this difficulty.
Third, the uncertainty sources are numerous: several tens of uncertain inputs have to be taken into account and screening techniques are required.
Finally, the model output variables of interest are several time-dependent curves (temperature, primary pressure, exchange coefficient) while well-known UASA are defined for scalar output variables.

In the following section, we present more deeply the industrial problem, the associated numerical model and the variables of interest.
In the third section, a functional screening technique, based on variance decomposition and generalized sensitivity indices, is described and applied.
In the fourth section, the functional metamodeling method is detailed.
Finally, a conclusion summarizes our results and gives some perspectives for this work.

\section{THE THERMAL-HYDRAULIC TRANSIENT SIMULATION}

In the generic methodology of uncertainty treatment (de~Rocquigny et al. \cite{derdev08}), the first step is devoted to the specification of the problem.
In particular, the objectives of the studies, the computer code which will be used, the studied scenario, the input uncertain variables and the output variables of interest have to be defined.

\subsection{The industrial problem}

\begin{figure}[!ht]
$$\psfig{figure=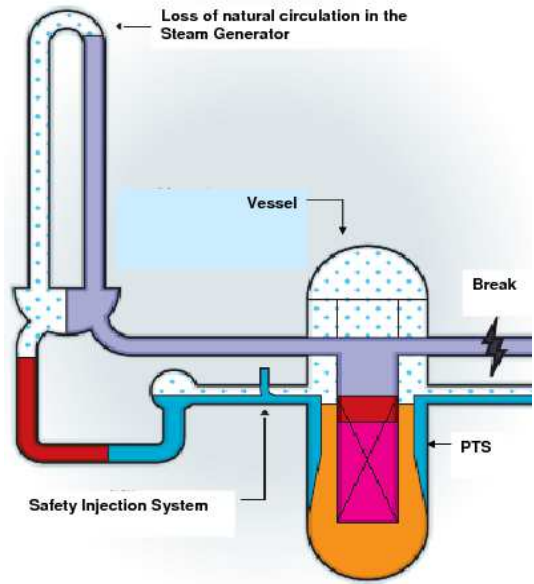,height=8cm,width=7cm}$$

\vspace{-0.4cm}
\caption{Scheme of a break in the hot leg of the primary circuit leading to a Pressurized Thermal Shock (PTS). Adapted from an IRSN document.}\label{fig:cuve}
\end{figure}

Evaluating nuclear reactor pressure vessel (NRPV) performance during transient
and accidental conditions is a major issue required by the regulatory authorities for the safety demonstration of the nuclear power plant.
Indeed, during the normal operation of a nuclear power plant, the NRPV walls are exposed to neutron radiation, resulting in localized embrittlement of the vessel steel and weld materials in the area of the reactor core.
Therefore, the knowledge of the behaviour of the pressure vessel subjected to highly hypothetic accidental conditions, as a Pressurized Thermal Shock (PTS) transient, is one of the most important inputs in the nuclear power plant lifetime program (see Fig. \ref{fig:cuve}). 

In our industrial context, the transient critical for the PTS consists of a Small Break LOCA (Loss of Coolant Accident), located at a hot leg of the studied PWR (see Fig. \ref{fig:cuve}). 
Consecutive water loss and depressurization lead to the initiation of the safety injection system. The cold water injected by this system may enter in contact with the hot walls of the vessel and may provoke a sudden temperature decrease (thermal shock) associated with severe related thermal stresses. The temperature decrease depends on the way the injected cold water mixes with the hot fluid and on the level of heat transfer between the more or less mixed cold water and the vessel. The safety issue corresponding to the PTS is the vessel rupture. Fluid temperature close to the vessel walls and heat transfer coefficient are, with the pressure inside the vessel, the three determinant variables in the prediction of PTS consequences. These three curves are time dependent and calculated by a thermal-hydraulic code. As many of the input variables of the thermal-hydraulic code are uncertain, it will be a challenging task to evaluate the uncertainties on these three time dependent responses. In the following of the paper, we will focus on only one of the responses of interest: the fluid temperature (named $T_{\mbox{bas}}$).

\subsection{The CATHARE2 model}

CATHARE is the French best-estimate code for thermal-hydraulic nuclear reactor safety studies. It is the result of a more than thirty years joint development effort by CEA (the research institution), EDF (the utility), AREVA (the vendor) and IRSN (the safety authority). 
The code is able to model any kind of experimental facility or PWR (Pressurized Water Reactors, the present Western type reactors), or is usable for other reactors (already built such as RBMK reactors, or still underway such as the future reactors of Generation IV). It is also used as a plant analyzer, in a full scope training simulator providing real time calculations for the plant operators.

The code has a modular structure, which allows all the quoted above modelings. The main module is the 1-D module, but there is also $0$-D and $3$-D volumes and all the modules can be connected to walls or heat exchangers. Many sub-modules are available to calculate for example the neutronics, the fuel thermo-mechanics, pump characteristics, accumulators, etc.  All modules use the two-fluid model to describe steam-water flows. Four non-condensable gases may also be transported. The thermal and mechanical non-equilibrium are described, as well as all kinds of two-phase flow patterns. 

A stringent process of validation is performed. First of all, quite analytical experiments are used in order to develop and qualify the physical models of the code. After that, a verification process on integral effect tests is carried out to verify the overall code performance. 

\subsection{Definition of uncertainty sources}

This section describes how the list of the input variables and the quantification of their uncertainty via a probability density function have been performed.
First of all, two existing lists were used as a basis. The first list is the one established by CEA for the BEMUSE benchmark (de Cr\'ecy et al. \cite{decbaz08}), devoted to the Large Break LOCA. The second list is aimed at defining the work program of determination by CEA of the uncertainties of the physical models of CATHARE, both for Small and Large Break LOCA. In both lists, the output variable of interest is the maximum clad temperature. The PTS is a transient which has some specificities, of course with respect to the Large Break LOCA, but also to the already studied Small Break LOCA, owing to the size and the position of the break. In addition to that, the output variables of interest in a PTS are related to the conditions of pressure, liquid temperature, flow rate and level in the cold leg and not to the maximum clad temperature. Consequently, a rather detailed analysis of the PTS scenario has also been performed by experts and has led to the definition of $p=31$ input variables, which are of two types.

The first type of variables is those related to the system behavior of the reactor. They are:
\begin{itemize}
\item	The initial conditions  (e.g. secondary circuit pressure);
\item	The boundary conditions (e.g. the residual power);
\item	The flow rate at the break (all the parameters involved in the critical flow rate prediction, e.g. liquid-interface heat transfer by flashing);
\item	The conditions upstream from the break (e.g. interfacial friction in stratified flow in the hot leg); 
\item	Interfacial friction (e.g. in the core);
\item	Bypass phenomenon in the vessel upper head.
\end{itemize}

The second type of variables is those related to the safety injections and accumulators. They are:
\begin{itemize}
\item	The temperature and the flow rate of the injections;
\item	The features of the accumulators (e.g. initial pressure);
\item	The condensation  phenomena (e.g. liquid-interface heat transfer in stratified flow downstream from the injections);
\item	The thermal stratification in the cold leg. 
\end{itemize}

The quantification of the uncertainty of the variables is performed for the physical models either by using a statistical method of data analysis, the name of which is CIRCE (giving normal or log-normal distribution types), or by expert judgment (giving also normal or log-normal distributions). For the initial and boundary conditions, reactor data are generally  used, giving uniform or normal laws. CIRCE (de Cr\'ecy \cite{dec01}) can be used only if there are experiments with a sufficient numbers of measures, which is the case of a rather large number of variables of both existing lists: BEMUSE and CATHARE. By lack of time, CIRCE was not used for the ``new'' parameters, specific to the PTS. A CIRCE study is nevertheless underway for two parameters related to the condensation effects resulting from the injection of cold water in the cold leg. 

At the end of this step, $31$ random input variables have been defined with their probability density function.
The following step consists in performing a sensitivity analysis in order to retain only the most influent inputs on the output variables of interest.

\section{SCREENING WITH GENERALIZED SENSITIVITY INDICES}\label{sec:screening}

In order to reduce the initial list of $31$ variables for the construction of metamodels and the propagation of uncertainties, it is necessary to implement a method that gives the sensitivity of the input variables on the responses (i.e. the output variables of interest) obtained with the code CATHARE2.
For a scalar response, this approach is called the screening and numerous methods can be applied (Saltelli et al. \cite{salcha00}, Dean \& Lewis \cite{dealew06}).
In this paper, the model response is a time curve and these classical approaches are not applicable.

A first sensitivity analysis, using a One-At-a-Time (OAT) method has provided a first idea of the influence of the input variables. 
OAT is a coarse sensitivity analysis process where each input value is moved once (for example from a minimal value to a maximal value), and its effect on the output value is analyzed.
In our case, $p+1=32$ simulation runs have been performed; then a subjective visual analysis between the $32$ output curves has allowed to detect $14$ influent inputs.
For several reasons, OAT is judged as a bad practice for sensitivity analysis (see Annoni \& Saltelli \cite{salann10}), but engineers still use this approach in order to verify their model.
In the following, we use a more rigorous method for the sensitivity analysis of inputs on an output curve.

\subsection{Generalized sensitivity indices method}

For a discrete time response, a sensitivity analysis could be performed separately at each calculation time step. Indeed in this case, the studied output is a scalar and the classical global sensitivity techniques are applicable. However, this technique has the disadvantage of introducing a high level of redundancy because of the strong correlations between responses from one time step to the next one. It may also miss important dynamic features of the response.  

Recently, Lamboni et al. \cite{lammak09,lammon10} have developed the method of generalized sensitivity indices ($GSI$), which summarizes directly the effect of each input variable on the time series output of a computer code. This method can be divided in four steps: 
\begin{enumerate}
\item 
$n$ simulations are carried out with the computer code (CATHARE2 in our case) by varying the $p$ input variables $(X^1, ..., X^p)$ randomly or according to a factorial experimental design. This gives $n$ time responses.
The set of outputs can be represented in the form of a $n \times T$ matrix, where $T$ is the number of time sampling points (in our case $T = 512$):
\begin{equation}
\mathbb{Y} = 
\begin{pmatrix} 
Y_1^1 & \cdots & Y_1^T \\
\vdots & \ddots & \vdots \\
Y_n^1 & \cdots & Y_n^T
\end{pmatrix} \;.
\end{equation}
Each column $Y^t$ in $\mathbb{Y}$ represents the simulated values of the output variable at a given time, while each row of $\mathbb{Y}$ is a transient simulation obtained for a given set of values of the input variables. 
\item
The principal components analysis (PCA) is used to decompose the whole variability (or total inertia) in $\mathbb{Y}$. The total inertia is defined as $I(\mathbb{Y}_c)=\mbox{trace}(^t\mathbb{Y}_c \mathbb{Y}_c)$, where $\mathbb{Y}_c$ is the matrix $\mathbb{Y}$ with each column centered around its mean. Thus $^t\mathbb{Y}_c \mathbb{Y}_c$  is the variance-covariance matrix of the columns of $\mathbb{Y}$. The PCA decomposition is based on the eigenvalues and eigenvectors of $^t\mathbb{Y}_c\mathbb{Y}_c$. Let $\lambda_1,\ldots,\lambda_T$ denote the obtained eigenvalues in decreasing order and let  $\mathbb{L}$ denotes the $T \times T$ matrix of eigenvectors, where each column $I_k$ is the eigenvector associated to $\lambda_k$. The $N \times T$ matrix $\mathbb{H}$ of principal components scores is obtained by: $\mathbb{H}=\mathbb{Y}_c \mathbb{L}$.  The columns of $\mathbb{H}$, $h_k$ for $k = 1,2,\ldots,T$, are mutually orthogonal linear combinations of the $\mathbb{Y}_c$ columns. By construction, $\mathbb{H}$ has the same total inertia as $\mathbb{Y}_c$, but its inertia is mainly concentrated in the first principal components scores.
\item
The third step in the method involves performing an ANOVA (ANalysis Of Variance) type decomposition of each principal components $h_k$:
\begin{equation}\label{eq:decomp}
SS(h_k) = SS_{1,k} + \ldots + SS_{i,k} + \ldots + SS_{p,k} + SS_{12,k} + \ldots + SS_{ij,k} + \ldots + SS_{p-1 p,k} + \ldots + SS_{1..p,k}
\end{equation}
where $SS$ represents the sum of square of the effects.
 $SS_{i,k}$ is the main effect of input variable $i$ on the principal component $k$, $SS_{ij,k}$ is the effect of the interaction between the input variables $i$ and $j$ on the principal component $k$, etc.
Note that this decomposition is always possible if the factorial design is orthogonal and is unique if the factorial design is complete (Montgomery \cite{mon04}).
The sensitivity indices $SI$ on the principal component $h_k$ are defined by:
 \begin{equation}
 SI_W(h_k) = \frac{SS_{W,k}(h_k)}{SS(h_k)} 
 \end{equation}
with $W = 1,\ldots,p$ for the first order sensitivity indices (main effects), $W= \{1,2\}, \ldots, \{p-1,p\}$ for the second order effects, etc.
 We have $SS_{W,k}(h_k)=\|S_W h_k\|^2$ where $S_W$ denotes the orthogonal projection matrix on the subspace associated to $W$.
\item
Finally the $GSI$ are calculated by summing the $SI_W(h_k)$ of the principal components weighted by the inertia $I_k$ associated with the component $h_k$:
 \begin{equation}
 GSI_W = \sum_{k=1}^T \frac{I_k}{I} SI_W(h_k) \;.
 \end{equation}
These $GSI$ measure the contribution of the terms W (of first order, second order, etc.) to the total inertia of the time responses. For practical purposes, only the first $T_0$ principal components ($T_0 \ll T$)  are kept such that $\displaystyle \sum_{k=1}^{T_0} I_k =\frac{x}{100}I$ with $x$ a given percentage ($99\%$ for example). 
\end{enumerate}

The dynamic coefficient of determination $R^2_t$  evaluates the quality of the approximations (truncations of the principal components and of the ANOVA decomposition) directly on the original time series (matrix $\mathbb{Y}$):
\begin{equation}\label{eq:R2t}
R^2_t = \frac{\sum_{i=1}^n \left( \tilde{y}_i^t-\bar{y}^t\right)^2}{\sum_{i=1}^n \left( y_i^t-\bar{y}^t\right)^2} \;\mbox{ for } t=1,\ldots, T \;,
\end{equation}
where $\bar{y}^t=\sum_{i=1}^n y_i^t$ and $\tilde{y}^t$ are the columns of the approximated matrix of outputs $\tilde{\mathbb{Y}}_c=\displaystyle \sum_W S_W \tilde{\mathbb{H}} \tilde{\mathbb{L}}^t$.
The matrix $\tilde{\mathbb{H}}$ (resp. $\tilde{\mathbb{L}}$) contains the first $T_0$ columns of $\mathbb{H}$ (resp. $\mathbb{L}$) and the summation on $W$ is restricted to those factorial terms in the approximation.
 A value of $R^2_t$  close to $1$ indicates that most of the inertia of $\mathbb{Y}$ was captured while a low $R^2_t$  means that the obtained $GSI$ should be interpreted with caution.

The complete decomposition (\ref{eq:decomp}) is possible only in the case of a full factorial design (untractable in our case because a full factorial design at $2$ levels and $31$ variables would require $2^{31} \simeq 215 \times 10^6$ simulations). 
As a consequence, we use a fractional factorial design (Montgomery \cite{mon04}), which enables only studying a limited number of effects. 
We focus our study to a two levels fractional factorial design of resolution IV, which can estimate without bias the first order effects of $p=31$ input variables with $n=2^{p-q}=2^6=64$ simulations ($q=25$). The generated experimental design selects $64$ points among the $2^{31}$ possible points of the two levels full factorial design.

\subsection{Results for the response $T_{\mbox{bas}}$}

\begin{figure}[!ht]
$$\psfig{figure=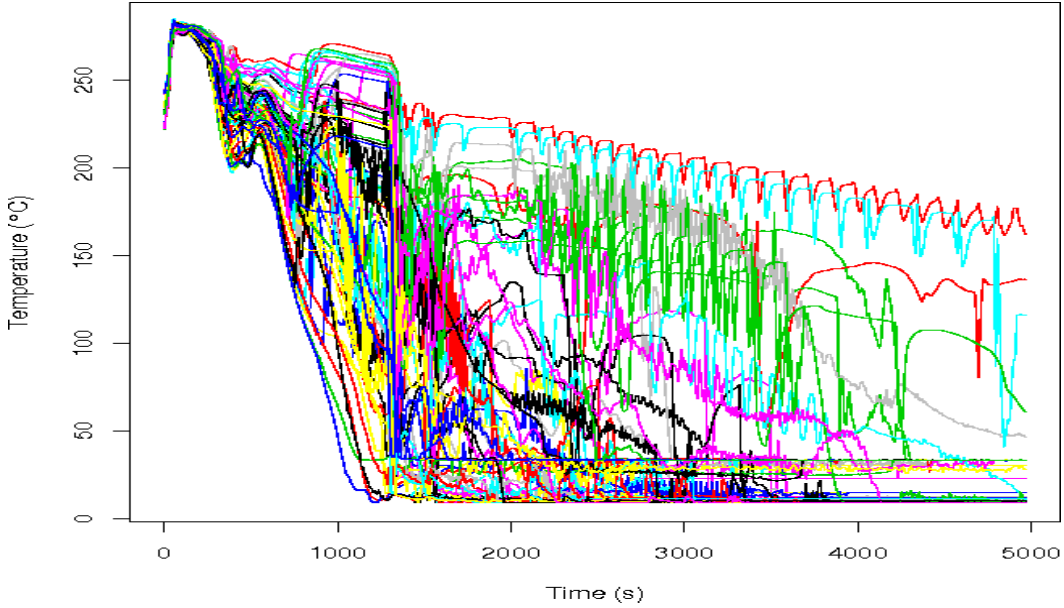,height=8cm,width=14cm}$$

\vspace{-0.4cm}
\caption{$64$ simulations of $T_{\mbox{bas}}$ obtained from the results of the CATHARE2 code.}\label{fig:Tbas}
\end{figure}

A two levels (corresponding to the selected minimal and maximal values of the variables) experimental design, and then with only $64$ points, was generated and the corresponding calculations performed with the CATHARE2 code. The results concern the response of interest $T_{\mbox{bas}}$, the temperature at the bottom of the cold leg as a function of time.  The $64$ times series are shown in the figure \ref{fig:Tbas}. $95\%$ of the inertia of these time series are captured by the nine first principal components (the three first principal components capture $82\%$ of the inertia).

\begin{figure}[!ht]
$$\psfig{figure=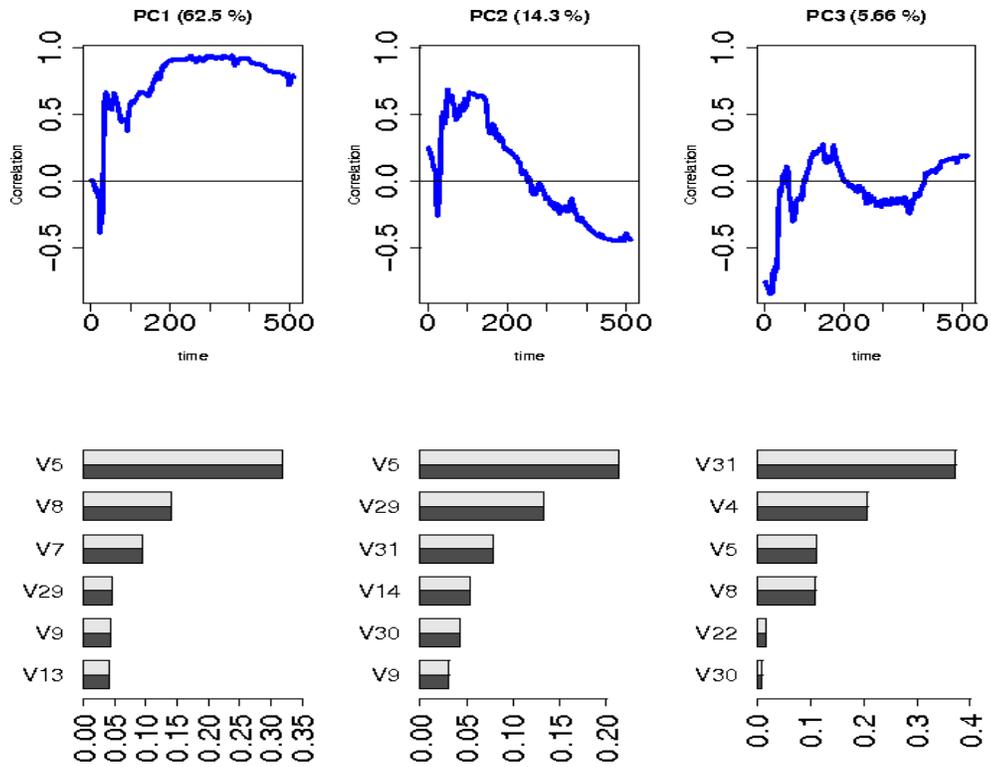,height=10cm,width=13cm}$$

\vspace{-0.4cm}
\caption{Up: correlation coefficients between the $3$ first principal components (PC) and the outputs $Y_t$; bottom: $6$ larger sensitivity indices $SI_W$ of the $3$ first PC.}\label{fig:3PC}
\end{figure}

Figure \ref{fig:3PC} shows at the top the correlation coefficients between the $3$ first principal components (PC) and the outputs $Y_t$ for $t =1,\ldots,T$ (with $T = 512$) and at the bottom the sensitivity indices $SI_W$ (only the six most important are shown) of the $3$ first PC. We notice that the first PC is positively correlated with the output except at the beginning of the transient, while the second has a better correlation at the beginning of the transient. The first two PC are more sensitive to the variable $5$ (residual power) and the third one to the variable $31$ (stratification rate). 

\begin{figure}[!ht]
$$\psfig{figure=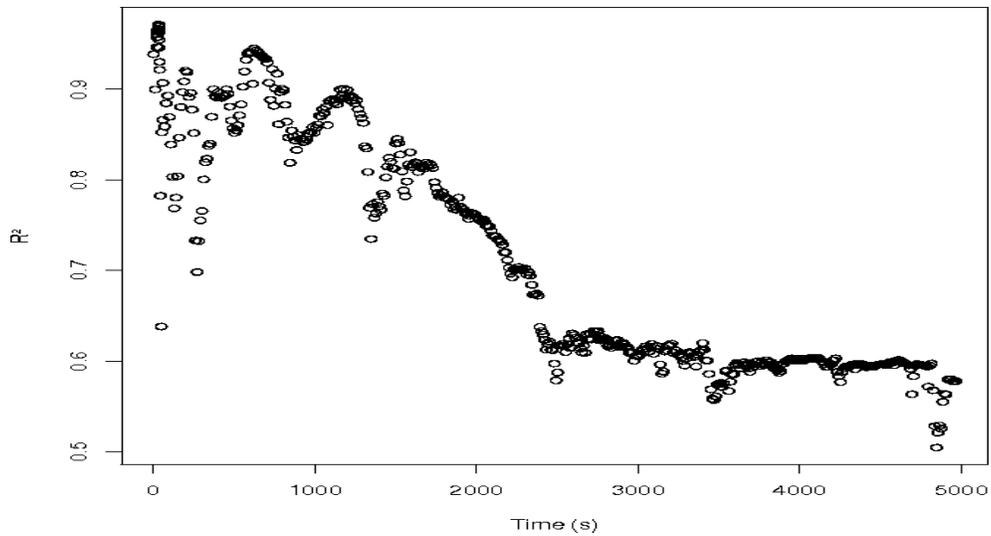,height=7cm,width=13cm}$$

\vspace{-0.4cm}
\caption{Dynamic coefficient of determination $R^2_t$ versus time.}\label{fig:dynR2}
\end{figure}

Figure \ref{fig:dynR2} shows the dynamic coefficient of determination versus time. The values of $R^2_t$, between $0.6$ and $1$, accredit a good confidence to the $GSI$ estimates. This confidence is better in the first part of the transient where $R^2_t$ are highest. This result is interesting because with regards to the thermal shock impact, we are more interested by the temperatures in the first part of the transient. We can therefore conclude that the $GSI$ estimates will be good indicators of the sensitivity of input variables to the thermal shock.

Table \ref{tab:screening} shows the $GSI$ of the variables, listed in decreasing order. The sum of the indices ($75\%$) represents the percentage of inertia (ie variability of the response) explained by the main effects. If we accept the variables having a $GSI$ greater than or equal to $1\%$, we obtain the $12$ variables listed in Table \ref{tab:screening}. When we compare these variables to the $13$ variables identified as important in the OAT analysis (see beginning of section \ref{sec:screening}), we find six variables in common: residual power ($V5$), steam generator emergency feedwater supply temperature ($V6$), injection system temperature ($V7$) and injection system flowrate ($V8$), accumulators temperature ($V9$) and nitrogen fraction in accumulators ($V13$). The other variables were not included in the OAT list: liquid-interface exchange in stratified flow ($V29$), rate of stratification ($V31$), singular pressure drop in dome bypass ($V14$), secondary circuit pressure ($V4$), liquid-interface exchange for all flow types ($V30$) and accumulators initial pressure ($V11$).
It shows that applying OAT can induce errors on the sensitivity analysis process and that the $GSI$ approach can produce useful results.

\begin{table}[!ht]
  \centering
\caption{First $GSI$ of the inputs (decreasing order).}\label{tab:screening}
\begin{tabular}{lll}
Input number & Input definition & $GSI (\%)$ \\
\hline\hline
$V5$	& Residual power & $25.45$ \\
\hline
$V8$	& Injection System flowrate	& $10.29$\\
\hline
$V7$	& Injection System temperature & $6.71$ \\
\hline
$V29$ & Liquid-interface exchange in stratified flow	& $6.03$\\
\hline
$V31$ & Rate of stratification	& $3.97$\\
\hline
$V9$	& Accumulators temperature & $3.35$\\
\hline
$V13$ & Nitrogen fraction in accumulators & $3.27$\\
\hline
$V6$ & Steam generator emergency feedwater supply temperature & $2.59$\\
\hline
$V14$ & Singular pressure drop in dome bypass & $2.51$\\
\hline
$V4$ & Secondary circuit pressure & $2.10$\\
\hline
$V30$ & Liquid-interface exchange for all flow types	& $1.74$\\
\hline
$V11$ & Initial pressure of accumulators & $1.03$
\end{tabular}
\end{table}
 
 In the following, we supress from this list of $12$ inputs one additional input (the nitrogen fraction in accumulators) and retain only $11$ input variables.
 This is due to a posteriori modeling considerations.

\section{THE FUNCTIONAL METAMODELING}

The section \ref{sec:screening} has introduced the $GSI$ which accuracy depends on the number of simulations.
Using a small number of simulations ($64$), this screening step has allowed us to find a restricted number of influent inputs on the output curves.
In this section, a methodology is proposed in order to build a metamodel of the computer code.
Not restricted by the cpu time, this metamodel will be able to be run as many times as necessary.
It will serve us for a quantitative sensitivity analysis and subsequent reliability studies.
Starting from $n$ input-output couples $(X_1,Y_1),\ldots,(X_n,Y_n) \in \R^p \times \R^T$, we first cluster the inputs-outputs into $K$ groups, then reduce the dimension $T$ to $d \ll T$ (for the reasons given in the paragraph 3.1), learn the relation between $x_i$ and reduced representations (using classical metamodeling techniques), and finally determine a reconstruction function from $\R^d$ to $\R^T$. 

\begin{figure}[!ht]
$$\psfig{figure=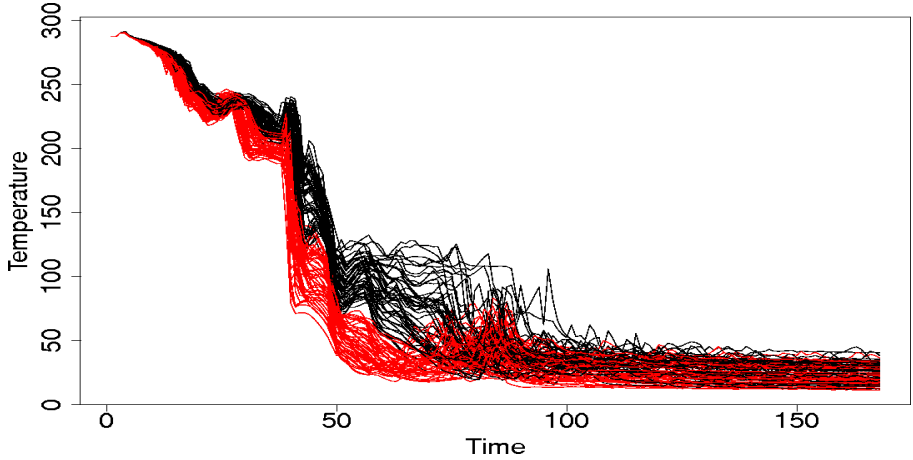,height=6cm,width=12cm}$$
\vspace{-0.8cm}
\caption{Temperature versus time (CATHARE2 code), two clusters.}\label{fig:catClust}
\end{figure}

The clustering step allows to distinguish several physical behaviors, which can be modeled in quite different ways. Our main ambition is to handle functional clusters of arbitrary shapes, and then reduce the dimension non-linearly to describe (almost) optimally any manifold.
Ascendant hierarchical clustering (with Ward linkage) will be the algorithm of choice, because it can detect classes of arbitrary shapes and is well suited to estimate the number of clusters. This last task is done by estimating prediction accuracy using subsampling, like in Roth et al. \cite{rotlan02}. We use the $L_2$ distance between curves to embed them into a $k$-nearest-neighbors graph, and then estimate the commute-time euclidian distances in this graph (Yen et al. \cite{yenvan05}). This preprocessing step eases the clustering process by increasing long distances and decreasing short ones. The figure \ref{fig:catClust} shows an example of two well separated clusters produced by the algorithm.

\subsection{Methods for dimensionality reduction}

Assume now that a group of $m < n$ inputs-outputs $(X_i, Y_i)$ has to be analyzed for dimension reduction (the indexation goes from 1 to $m$, but represents a subset of $\{1,\ldots,n\}$). The classical method for the dimension reduction is a (functional) principal components analysis (Ramsay \& Silverman \cite{ramsil05}), which is perfectly suited for dimensionality reduction of linear subspaces embedded into $\R^T$. However, in some cases this method is inappropriate; for example a simple sphere in $\R^3$ cannot be represented in two dimensions using principal components. Functions in $\R^T$ could have an even more complex structure ; that is why we choose to use a manifold learning technique, namely Riemannian Manifold Learning (RML) (Lin et al. \cite{linzha06}). This method aims at approximating Riemannian normal coordinates (see for example do~Carmo \cite{doc92}). It performs the best in our practical tests versus other nonlinear techniques.

The RML algorithm consists at the beginning in building a graph representing the data; all the further computations are done inside this graph. It is basically a $k$-nearest-neighbors graph, which is used to constrain the representations step by step. In short, RML ensures that the radial distances along the manifold, as well as the angles in each neighborhood are (approximately) preserved.

\subsubsection{Building of the initial graph}

The algorithm proceeds as follows:
\begin{enumerate}
\item{Build the $k$-nearest neighbor graph for some value of $k$ (around $\sqrt n$ where $n$ is the number of data points).}
\item{In each neighborhood, delete all edges not matching the visibility condition below.}
\item{Whenever two connected points are relatively too far apart, delete the corresponding edge.}
\end{enumerate}
The last step is rather heuristic, and could need to be adapted to the data. That is why we do not describe it there. 

\begin{figure}[!ht]
$$\psfig{figure=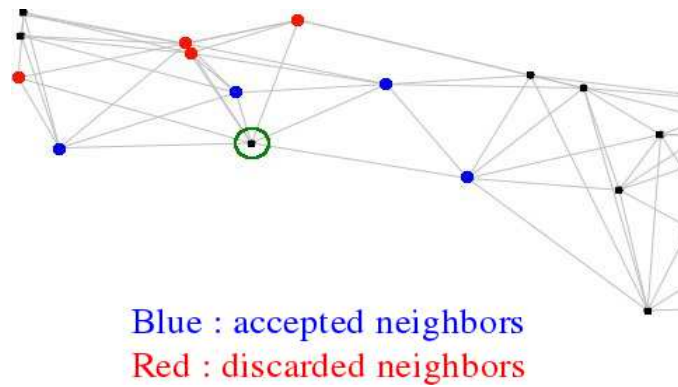,height=5cm,width=9cm}$$
\vspace{-0.8cm}
\caption{Visibility neighborhood of a point $Y$ (circle).}\label{fig:knnvois}
\end{figure}

For the second step, we define $VN$ as the ``visible neighbors'' set of $Y$.
We laso define the function $\text{visibility}(VN,i)$ which is true if every index $i_{\ell}$ in $VN$ satisfies $\widehat{Y Y_{i_{\ell}} Y_i} \leq \frac{\pi}{2}$, and which is false otherwise.
This means that the neighbors should not "eclipse" each other, which is a natural condition to avoid edge redundancy and thus better models the local geometry. 
Considering a neighborhood $Y_{i_1},\dots,Y_{i_k}$ of $Y$ in $\{Y_i, i=1,\dots,n\}$, the second step proceeds as follow:
\begin{enumerate}
\item{Set $VN$ to $Y_{i_1}$.}
\item{For $i$ from 2 to $k$:\\
\hspace*{1cm}if $\text{visibility}(VN,i)$, add $Y_i$ to $VN$.}
\item{Set the neighbors of $Y$ to $VN$.}
\end{enumerate}
The figure \ref{fig:knnvois} summarizes the neighbors selection process.

\subsubsection{Determining the starting point and initial coordinates}

Quite heuristically, $Y_0$ will stand for the point which minimizes the sum of graph distances from it to all other points (we use the Dijsktra algorithm to find it). This is a way to choose the ``center'' of the functional data, which is expected to be a good starting point. Now, to initialize the stepwise search for coordinates $Z_i \in \R^d$, we need to know the coordinates of a small neighborhood around $Y_0$. Practically, this is done by computing a principal component analysis on the selected small neighborhood (of size $k$ or less). This provides $d$ basis functions, which generate a $d$-dimensional space $Q_0$.

\begin{figure}[!ht]
$$\psfig{figure=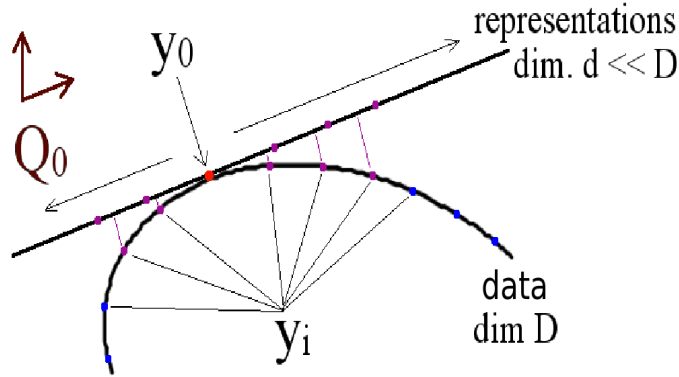,height=5cm,width=9cm}$$
\vspace{-0.8cm}
\caption{Local step in the RML algorithm.}\label{fig:rmlStep1}
\end{figure}

If $B_0$ is the matrix of the $d$ discretized basis functions in lines, then the coordinates $Z_i$ when $Y_i$ is a neighbor of $Y_0$ are initially determined as follows:
$$Z_i = B_0 (Y_i - Y_0) \, .$$
That is, $Z_i$ is the coordinates of the projection of $Y_i$ onto $Q_0$. A final normalization ensures that radial distances are preserved: $Z_i$ is replaced by $\frac{\|Y_i-Y_0\|}{\|Z_i\|} Z_i$ (with $Z_0 = 0$). 
This is illustrated on the figure \ref{fig:rmlStep1}.

\subsubsection{Solving for all coordinates by local optimization}

If $Y$ is too far from $Y_0$ for a valid tangent plane approximation, the basic idea is to preserve local angles with the neighbors of a predecessor on a shortest path from $Y_0$ to $Y$. Indeed, this results in an (approximately) locally conformal application. If we write $Y_p$ for a predecessor of $Y$, and $Y_p^{(1)},\ldots,Y_p^{(q)}$ its neighbors with known reduced coordinates $Z_p^{(i)}$ (some are surely known if we use a breadth-first walk in the graph), we have to ensure that $\widehat{Z Z_p Z_p^{(i)}} \simeq \widehat{Y Y_p Y_p^{(i)}}, i=1,\ldots,q$ under the distance constraint $\|Z-Z_p\| = \|Y-Y_p\|$. This can be transformed into a least-squares problem with quadratic constraints, leading to an equation solved numerically. The procedure is detailed below.

\begin{figure}[!ht]
$$\psfig{figure=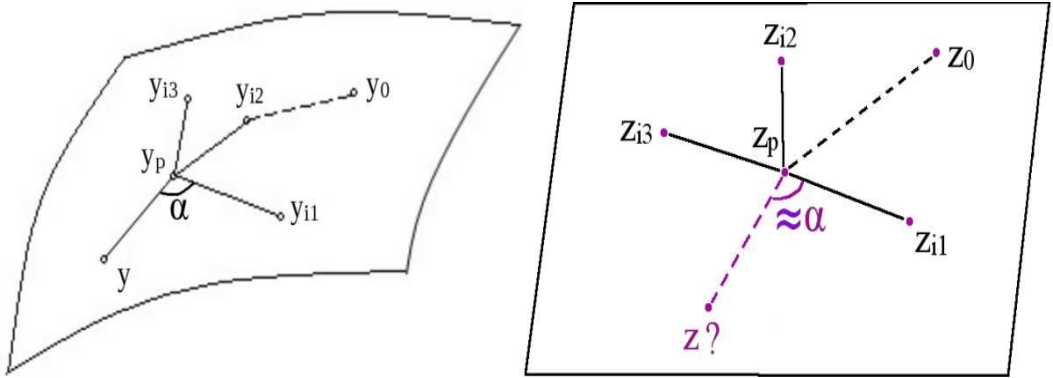,height=5cm,width=14cm}$$
\vspace{-0.8cm}
\caption{Global step in the RML algorithm. $Y_p$ (resp. $Z_p$) is a predecessor of $Y$ (resp. $Z$) and $(Y_{i1},Y_{i2},Y_{i3})$ (resp. $(Z_{i1},Z_{i2},Z_{i3})$) are the neighbors of $Y_p$.}\label{fig:rmlStep2}
\end{figure}

Given three vectors $\alpha$, $\beta$ and $\gamma$, 
$\cos \widehat{\alpha \beta \gamma}$ can be written $\frac{\langle \alpha-\beta, \gamma-\beta \rangle}{\|\alpha-\beta\| . \|\gamma-\beta\|}$. Thus, the optimization problem can be transformed in
$$\forall i=1,\dots,n \, , \, \frac{\langle Z-Z_p, Z_p^{(i)}-Z_p \rangle}{\|Z-Z_p\|. \|Z_p^{(i)}-Z_p\|} - \frac{\langle Y-Y_p, Y_p^{(i)}-Y_p \rangle}{\|Y-Y_p\| . \|Y_p^{(i)}-Y_p\|} \simeq 0 \, ,$$
under the constraint $\|Z-Z_p\| = \|Y-Y_p\|$. The "$\simeq$" sign above has to be specified ; we choose to minimize the euclidian norm of the vector of the angles differences. By writing $A_Z$ for the matrix of the vectors $\frac{Z_p-Z_p^{(i)}}{\|Z_p-Z_p^{(i)}\|}$ in lines, and $B_Y$ for vector of dot products $\frac{\langle Y-Y_p, Y_p^{(i)}-Y_p \rangle}{\|Y_p^{(i)}-Y_p\|}$ (fully known), the optimization problem can be stated compactly as follows:
$$\mbox{Minimize } \|A_Z (Z - Z_p) - B_Y\| \, ,$$
under the constraint $\|Z - Z_p\| = \|Y - Y_p\|$. We can then write the Lagrangian of this optimization problem, to express $Z - Z_p$ in function of all known variables and the Lagrange multiplier $\lambda$. This last parameter can be estimated at low cost by solving an equation numerically. 
The figure \ref{fig:rmlStep2} illustrates this process.

The reconstruction function is then easily done. We start from the training sample of functional outputs $Y_i \in \R^T$ and corresponding embeddings $Z_i \in \R^d$. 
Let us write $Z_{i_1},\ldots,Z_{i_k}$ the $k$ nearest neighbors of $Z$ in the $d$-dimensional space. First, a local functional PCA is achieved around $Y_{i_1},\ldots,Y_{i_k}$ to replace the curves by their $d$-dimensional principal components scores $Y'_{i_j}$. Then we solve the problem stated above with reversed roles, before re-expanding reduced coordinates $Y'$ onto the orthonormal basis of (functional) principal components.

\subsection{Application to the PTS}

\begin{figure}[!ht]
$$\psfig{figure=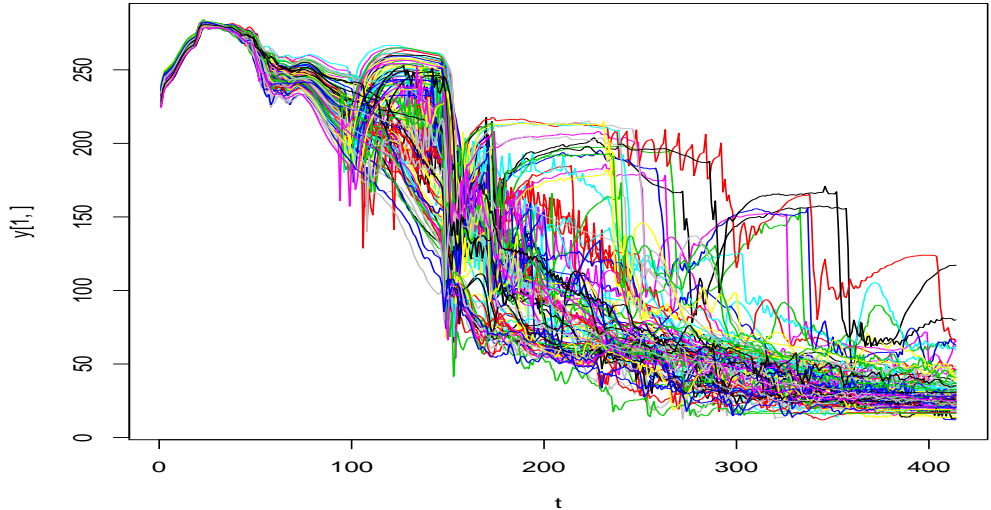,height=8cm,width=14cm}$$
\vspace{-1.1cm}
\caption{$100$ curves amonst the $600$ output curves.}\label{fig:cathare2}
\end{figure}

A special low discrepancy design of size $n=600$ with the $11$ influent input variables has been built.
This design is an optimized Latin Hypercube Sample with a minimal wrap-around discrepancy, which has been proved to be especially efficient for the metamodel fitting process (Iooss et al. \cite{ioobou10}).
Therefore, $600$ simulation runs have been performed in order to obtain the $600$ output curves.
Each curve is sampled on $T=414$ points (see Fig. \ref{fig:cathare2} for an example of $100$ curves).

We compare three kinds of metamodels: functional PCA, (functional) RML and a direct estimate through a local weighted mean, which does not use any dimensionality reduction. This last model, written F$k$NN for "functional $k$ nearest neighbors", predicts an output $Y \in \R^T$ from a new input $X \in \R^p$ as follows:
\begin{enumerate}
\item{Determine the $k$ nearest neighbors of $X$ among the set $\{X_1,\dots,X_n\}$ ; write their indices $i_1,\dots,i_k$.}
\item{Compute ${\displaystyle \hat Y = \frac{1}{C} \sum_{j=1}^{k} e^{\frac{-\|X-X_{i_j}\|^2}{\sigma_{i_j}^2} Y_{i_j}}}$, where $C = \sum_{j=1}^{k} e^{\frac{-\|X-X_{i_j}\|^2}{\sigma_{i_j}^2}}$.}
\end{enumerate}
The main parameter, $k$, is estimated by cross-validation. The local parameters $\sigma_i$ are computed by maximizing the local variations among gaussian similarities, according to the formula
$$\sigma_i = \frac{\|X_i-X_{i_k}\|^2 - \|X_i-X_{i_1}\|^2}{\log \|X_i-X_{i_k}\|^2 - \log \|X_i-X_{i_1}\|^2} \, .$$
For the other models (functional PCA and RML), we use the Projection Pursuit Regression method (see for example Friedman \& Stuetzle \cite{fristu81}) as a metamodel to learn the reduced coefficients. 

\begin{figure}[!ht]
$$\psfig{figure=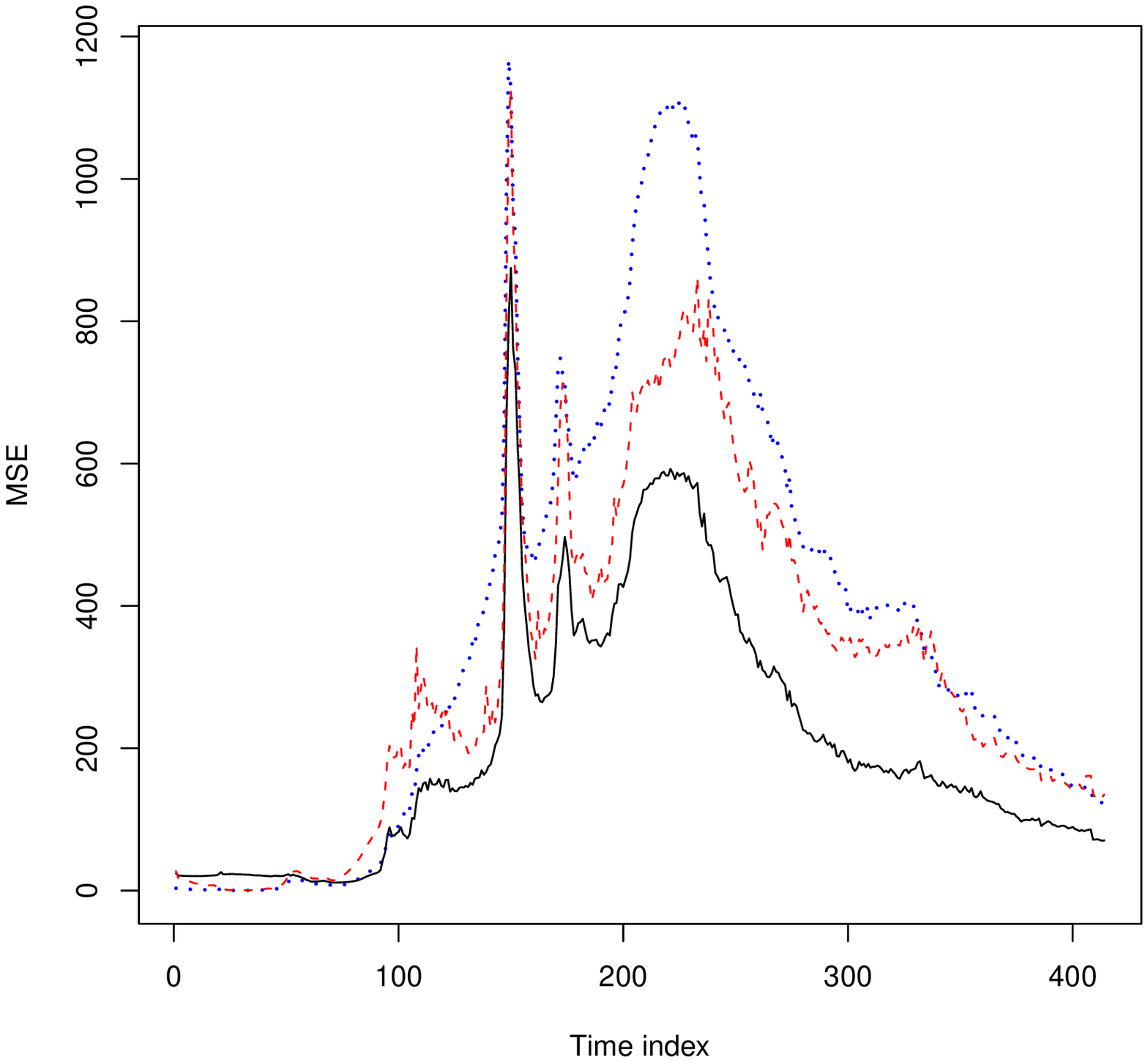,height=6cm,width=8cm}
\psfig{figure=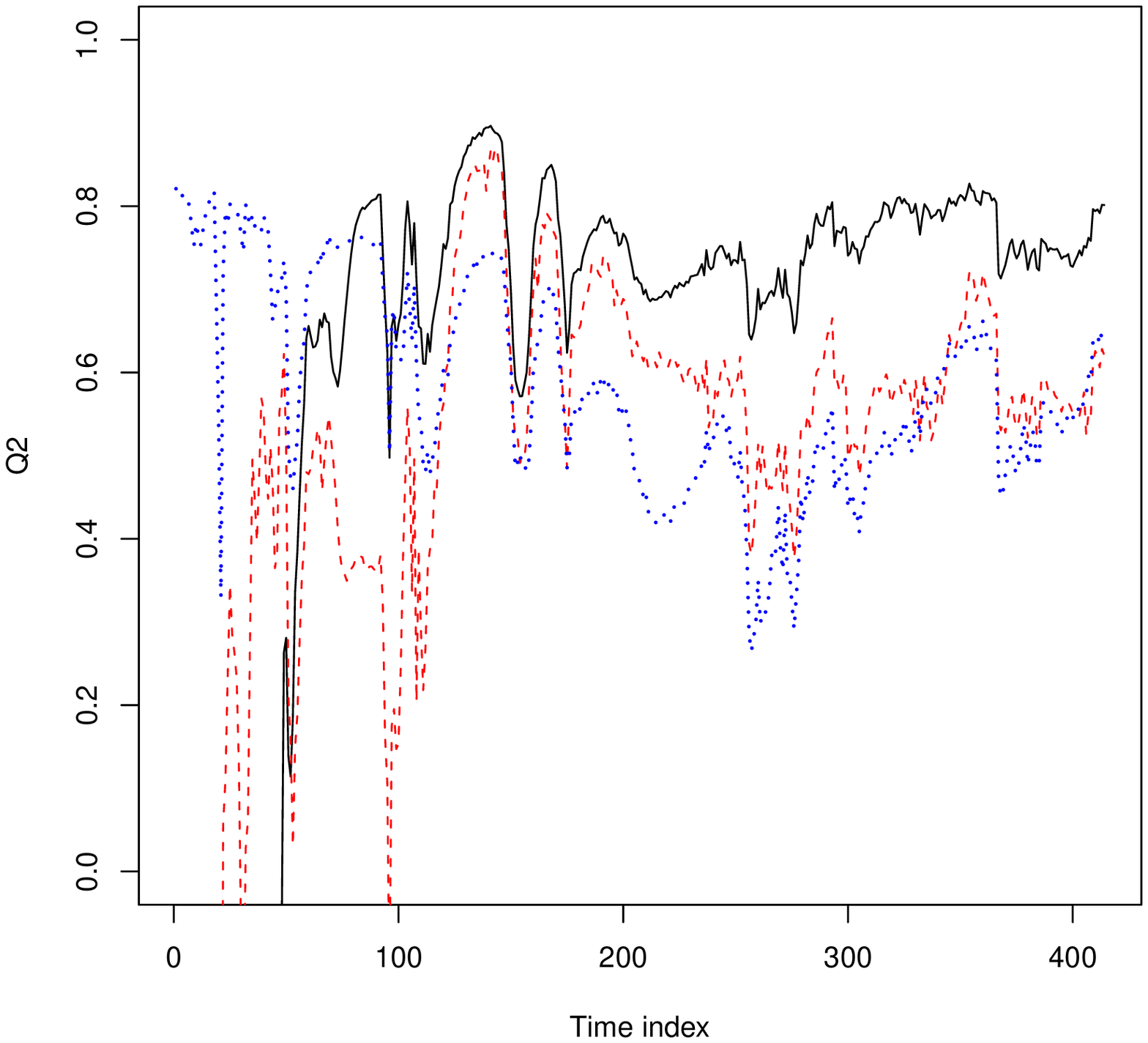,height=6cm,width=8cm}$$
\vspace{-1.1cm}
\caption{Mean square error $MSE$ (left) and predictivity coefficient $Q_2$ (right) error curves (solid line: PCA, dashed line: RML, dotted line: F$k$NN).}\label{fig:MSEQ2}
\end{figure}

Figure \ref{fig:MSEQ2} shows the mean square error ($MSE$) and the predictivity coefficient $Q_2$ (coefficient of determination computed on test samples) pointwise errors.
These measures have been computed using a $10$-fold cross-validation technique (dividing the $600$ simulation sets in $10$ simulation sets).
A $Q_2$ small or negative corresponds to a poor model, while a $Q_2$ close to $1$ corresponds to a good model. 
The small $Q_2$ values for RML and PCA at the beginning of the transient (time index smaller than $100$) are non significative because they are due to the small variability of the curves, turning the output variance to be close to $0$ (which is the denominator in the $Q_2$ formula).
The corresponding small values of the $MSE$ confirm this behaviour. 
For Nadaraya-Watson curve, this effect is not present because this estimator computes a local mean and the $MSE$ at the beginning of the curves are extremely small.
For the other parts of the curve (time index larger than $100$), the F$k$NN estimate is clearly insufficient compared to the PCA and RML estimates.
Globally, the PCA method performs best: the predictivity coefficients are larger than $70\%$ (for time index larger than $100$) which is clearly a good result with regard to the complexity and chaotic aspect of the CATHARE2 output curves.
%

\begin{figure}[!ht]
$$\psfig{figure=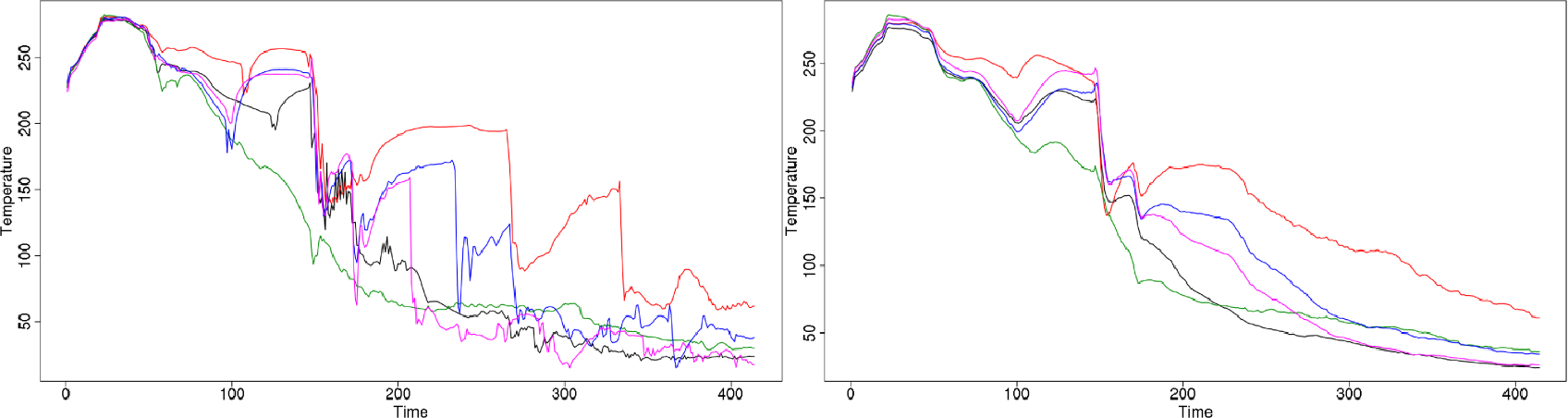,height=5cm,width=15cm}$$
$$\psfig{figure=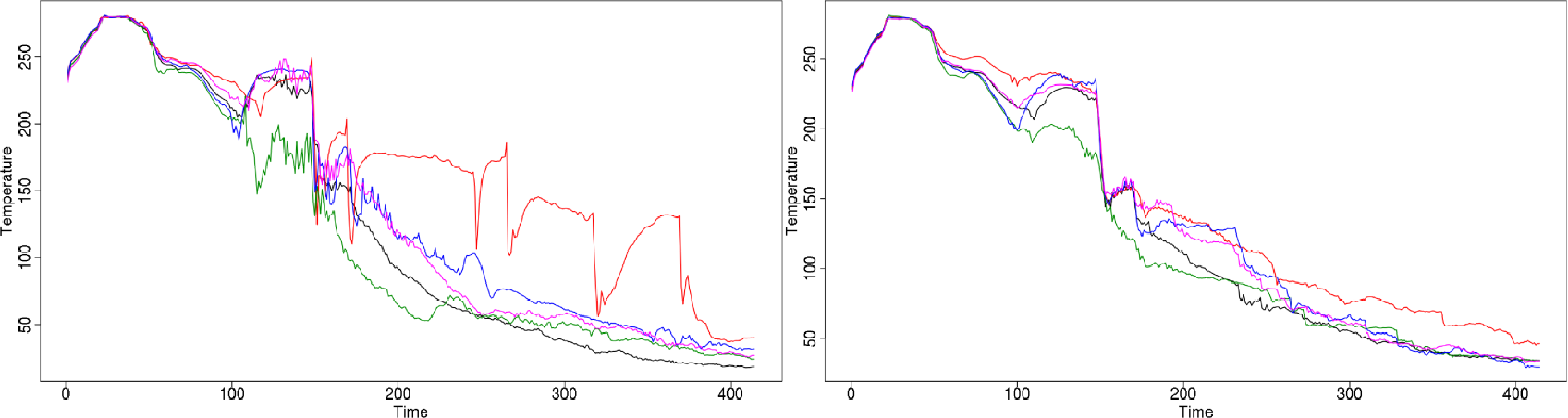,height=5cm,width=15cm}$$

\vspace{-0.5cm}
\caption{Up left: $5$ original (CATHARE2 simulation) curves. Up right: corresponding predicted curves with functional PCA. Bottom left: corresponding predicted curves with RML. Bottom right: corresponding predicted curves with F$k$NN method.}\label{fig:catSample}
\end{figure}

To visualize more accurately the performances of the different methods, five metamodel-based predicted curves are given in Figure \ref{fig:catSample}.
Several curves are better approximated by the RML method, but several others are really badly predicted by RML comparing to the PCA predictions.
Indeed, RML reproduces the form of the curves irregularities, but these irregularities are sometimes misplaced. 
We have seen previously that PCA gives better results than RML in terms of $MSE$ and $Q_2$ criteria.
In fact, these $L_2$ criteria are not the most adapted ones to measure the adequacy of the predicted curves when we are interested by reproducing the irregularities and the jumps of the curves.
As a conclusion, linear (PCA) and nonlinear (RML) techniques seem complementary. 
We are still working on improvements of the RML technique in order to improve its $L_2$ performance.

\section{CONCLUSION}

In order to improve the UASA of the PTS analysis using probabilistic methods, a methodology of metamodel fitting on the output curves of a thermal-hydraulic computer code has been proposed in this paper.
Compared to other works on this subject which deal with scalar outputs (de~Rocquigny et al. \cite{derdev08}, Saltelli et al. \cite{salrat08}), considering the whole output curve renders the UASA more difficult.
We have shown that the GSI approach of Lamboni et al \cite{lammak09} is particularly adapted to a screening process by using a fractional factorial design technique.
For the metamodeling process, due to the extremely complex behaviour of the thermal-hydraulic output curves, a new methodology has been developed in several steps: curves clustering, dimensionality reduction, metamodel fitting and curves reconstruction.
Even if the functional PCA and RML methods of dimensionality reduction have provided complementary and promising results, they need further improvements to precisely capture the irregularities of the output curves.
One first possibility would be to mix these approaches in order to keep the best properties of each one.


The final goal of a functional metamodel development is to be able to take into account the thermal-hydraulic uncertainties when performing a reliability analysis of the PTS scenario.
Indeed, such analysis needs several thousands of model runs in order to compute the vessel failure probability subject to this scenario.
Using a best-estimate thermal-hydraulic computer code as CATHARE2 (required by such safety analyses), this large number of runs remains intractable.
This paper puts a first shoulder to the wheel by proposing to approximate the thermal-hydraulic computer code by a metamodel.
Of course, as our metamodels have been built in the goal of a global approximation, they cannot be directly used for performing reliability analysis.
In the latter case the metamodel shall be accurate only in the region of the space of input variables where the performance function is close to zero and not everywhere.
For this issue, several methods have been developed in the context of scalar metamodel (Bect et al. \cite{becgin11}), but none in the functional metamodel case.
Future works will try to estimate and control the uncertainty of the functional metamodel, then to propagate it inside the reliability problem.

Finally, our methodology can be applied to many other industrial fields and UASA applied problems, as in automotive, aerospace engineering, environmental sciences, etc., where the computer codes are used for prediction of temporal phenomena.

\section{ACKNOWLEDGMENTS}

Part of this work has been backed by French National Research Agency (ANR) through COSINUS program (project COSTA BRAVA n$^o$ANR-09-COSI-015). 
The application concerning the PTS has been performed within the framework of a collaboration supported by CEA/DISN, EDF and AREVA-NP.
We thank G\'erard Biau for helful discussions.
All the statistical parts of this work have been performed within the R environment.
We are grateful to Herv\'e Monod and Matieyendou Lamboni who have provided the ``planor'' and``multisensi'' R packages to compute respectively the fractional factorial design and the generalized sensitivity indices.
We have also used the ``modelcf'' R package (author: B. Auder) which contains functional metamodel functions.


\singlespacing
\bibliographystyle{plain}

\begin{thebibliography}{}

\end{thebibliography}


\begin{thebibliography}{10}

\bibitem{bayber07}
M.J. Bayarri, J.O. Berger, J.~Cafeo, G.~Garcia-Donato, F.~Liu, J.~Palomo, R.J.
  Parthasarathy, R.~Paulo, J.~Sacks, and D.~Walsh.
\newblock Computer model validation with functional output.
\newblock {\em The Annals of Statistics}, 35:1874--1906, 2007.

\bibitem{becgin11}
J.~Bect, D.~Ginsbourger, L.~Li, V.~Picheny, and E.~Vazquez.
\newblock Sequential design of computer experiments for the estimation of a
  probability of failure.
\newblock {\em Statistics and Computing, in press}, 2012.

\bibitem{cammck06}
K.~Campbell, M.D. McKay, and B.J. Williams.
\newblock Sensitivity analysis when model ouputs are functions.
\newblock {\em Reliability Engineering and System Safety}, 91:1468--1472, 2006.

\bibitem{dec01}
A.~de~Cr\'ecy.
\newblock Determination of the uncertainties of the constitutive relationships
  of the {CATHARE} 2 code.
\newblock In {\em M\&C 2001}, Salt Lake City, Utah, USA, september 2001.

\bibitem{decbaz08}
A.~{de Cr\'ecy}, P.~Bazin, H.~Glaeser, T.~Skorek, J.~Joucla, P.~Probst,
  K.~Fujioka, B.D. Chung, D.Y. Oh, M.~Kyncl, R.~Pernica, J.~Macek, R.~Meca,
  R.~Macian, F.~D'Auria, A.~Petruzzi, L.~Batet, M.~Perez, and F.~Reventos.
\newblock Uncertainty and sensitivity analysis of the {LOFT L2-5} test:
  {R}esults of the {BEMUSE} programme.
\newblock {\em Nuclear Engineering and Design}, 12:3561--3578, 2008.

\bibitem{der06a}
E.~{de Rocquigny}.
\newblock La ma\^{\i}trise des incertitudes dans un contexte industriel -
  1\`ere partie~: une approche m\'ethodologique globale bas\'ee sur des
  exemples.
\newblock {\em Journal de la Soci\'et\'e Fran\c{c}aise de Statistique},
  147(3):33--71, 2006.

\bibitem{derdev08}
E.~{de Rocquigny}, N.~Devictor, and S.~Tarantola, editors.
\newblock {\em Uncertainty in industrial practice}.
\newblock Wiley, 2008.

\bibitem{dealew06}
A.~Dean and S.~Lewis, editors.
\newblock {\em Screening - Methods for experimentation in industry, drug
  discovery and genetics}.
\newblock Springer, 2006.

\bibitem{doc92}
M.~P. do~Carmo.
\newblock {\em Riemannian Geometry}.
\newblock Birkhauser, 1992.

\bibitem{fanli06}
K-T. Fang, R.~Li, and A.~Sudjianto.
\newblock {\em Design and modeling for computer experiments}.
\newblock Chapman \& Hall/CRC, 2006.

\bibitem{fristu81}
J.~H. Friedman and W.~Stuetzle.
\newblock {Projection Pursuit Regression}.
\newblock {\em Journal of the American Statistical Association}, 76:817--823,
  1981.

\bibitem{heljoh06}
J.C. Helton, J.D. Johnson, C.J. Salaberry, and C.B. Storlie.
\newblock Survey of sampling-based methods for uncertainty and sensitivity
  analysis.
\newblock {\em Reliability Engineering and System Safety}, 91:1175--1209, 2006.

\bibitem{ioobou10}
B.~Iooss, L.~Boussouf, V.~Feuillard, and A.~Marrel.
\newblock Numerical studies of the metamodel fitting and validation processes.
\newblock {\em International Journal of Advances in Systems and Measurements},
  3:11--21, 2010.

\bibitem{ioovan06}
B.~Iooss, F.~{Van Dorpe}, and N.~Devictor.
\newblock Response surfaces and sensitivity analyses for an environmental model
  of dose calculations.
\newblock {\em Reliability Engineering and System Safety}, 91:1241--1251, 2006.

\bibitem{kle08}
J.P.C. Kleijnen.
\newblock {\em Design and analysis of simulation experiments}.
\newblock Springer, 2008.

\bibitem{lammak09}
M.~Lamboni, D.~Makowski, S.~Lehuger, B.~Gabrielle, and H.~Monod.
\newblock Multivariate global sensitivity analysis for dynamic crop models.
\newblock {\em Fields Crop Research}, 113:312--320, 2009.

\bibitem{lammon10}
M.~Lamboni, H.~Monod, and D.~Makowski.
\newblock Multivariate sensitivity analysis to measure global contribution of
  input factors in dynamic models.
\newblock {\em Reliability Engineering and System Safety}, 96:450--459, 2011.

\bibitem{linzha06}
T.~Lin, H.~Zha, and S.~U. Lee.
\newblock Riemannian manifold learning for nonlinear dimensionality reduction.
\newblock In {\em 9th European Conference on Computer Vision}, pages 44--55,
  Graz, Austria, May 2006. Springer.

\bibitem{madkre86}
H.O. Madsen, S.~Krenk, and N.C. Lind, editors.
\newblock {\em Methods of structural safety}.
\newblock Prentice Hall, 1986.

\bibitem{marpig05}
M.~Marqu\`es, J.F. Pignatel, P.~Saignes, F.~D'Auria, L.~Burgazzi, C.~M\"uller,
  R.~Bolado-Lavin, C.~Kirchsteiger, V.~La Lumia, and I.~Ivanov.
\newblock Methodology for the reliability evaluation of a passive system and
  its integration into a probabilistic safety assesment.
\newblock {\em Nuclear Engineering and Design}, 235:2612--2631, 2005.

\bibitem{marioo10}
A.~Marrel, B.~Iooss, M.~Jullien, B.~Laurent, and E.~Volkova.
\newblock Global sensitivity analysis for models with spatially dependent
  outputs.
\newblock {\em Environmetrics}, 22:383--397, 2011.

\bibitem{mon04}
D.C. Montgomery.
\newblock {\em Design and analysis of experiments}.
\newblock John Wiley \& Sons, 6th edition, 2004.

\bibitem{mungar10}
M.~Munoz-Zuniga, J.~Garnier, E.~Remy, and E.~de~Rocquigny.
\newblock Adaptive directional stratification for controlled estimation of the
  probability of a rare event.
\newblock {\em Reliability Engineering and System Safety}, submitted, 2010.

\bibitem{ramsil05}
J.~O. Ramsay and B.~W. Silverman.
\newblock {\em {Functional Data Analysis}}.
\newblock Springer, 2005.

\bibitem{rotlan02}
V.~Roth, T.~Lange, M.~Braun, and J.~Buhmann.
\newblock A resampling approach to cluster validation.
\newblock In {\em {International Conference on Computational Statistics}},
  volume~15, pages 123--128, Berlin, Germany, 2002. Springer.

\bibitem{sacwel89}
J.~Sacks, W.J. Welch, T.J. Mitchell, and H.P. Wynn.
\newblock Design and analysis of computer experiments.
\newblock {\em Statistical Science}, 4:409--435, 1989.

\bibitem{salann10}
A.~Saltelli and P.~Annoni.
\newblock How to avoid a perfunctory sensitivity analysis.
\newblock {\em Environmental Modelling and Software}, 25:1508--1517, 2010.

\bibitem{salcha00}
A.~Saltelli, K.~Chan, and E.M. Scott, editors.
\newblock {\em Sensitivity analysis}.
\newblock Wiley Series in Probability and Statistics. Wiley, 2000.

\bibitem{salrat08}
A.~Saltelli, M.~Ratto, T.~Andres, F.~Campolongo, J.~Cariboni, D.~Gatelli,
  M.~Salsana, and S.~Tarantola.
\newblock {\em Global sensitivity analysis - The primer}.
\newblock Wiley, 2008.

\bibitem{sob93}
I.M. Sobol.
\newblock Sensitivity estimates for non linear mathematical models.
\newblock {\em Mathematical Modelling and Computational Experiments},
  1:407--414, 1993.

\bibitem{volioo08}
E.~Volkova, B.~Iooss, and F.~{Van Dorpe}.
\newblock Global sensitivity analysis for a numerical model of radionuclide
  migration from the {RRC} "{K}urchatov {I}nstitute" radwaste disposal site.
\newblock {\em Stochastic Environmental Research and Risk Assesment},
  22:17--31, 2008.

\bibitem{yenvan05}
L.~Yen, D.~Vanvyve, F.~Wouters, F.~Fouss, M.~Verleysen, and M.~Saerens.
\newblock Clustering using a random-walk based distance measure.
\newblock In {\em {Symposium on Artificial Neural Networks}}, volume~13, pages
  317--324, Bruges, Belgium, 2005.

\end{thebibliography}

\end{document}